\documentclass[12pt]{amsart}

\usepackage{geometry}        
\usepackage{memhfixc}  
\usepackage{graphicx}
\usepackage{amssymb}
\usepackage{epstopdf}
\usepackage{amscd}
\usepackage{tikz}
\usepackage[colorlinks=true, pdfstartview=FitV, linkcolor=blue, 
            citecolor=blue, urlcolor=blue]{hyperref}

\usepackage{memhfixc}  
\usepackage{tikz}
\usepackage{pdfsync}  
\usepackage{amsmath,amsthm,indentfirst,marvosym}
\usepackage{amssymb}
\usepackage{amsfonts}
\theoremstyle{plain}

\usepackage[utf8]{inputenc} 
\usepackage{textcomp} 
\usepackage{graphicx}  
\usepackage{flafter}  
\usepackage{bm}  
\usepackage[T1]{fontenc}

\usepackage{tikz}
\usepackage{color}
\usepackage{memhfixc}  
\usepackage{pdfsync}  
\usepackage{amsmath,amsthm,indentfirst,marvosym}
\usepackage{amssymb}
\usepackage{amsfonts}
\theoremstyle{plain}
\newtheorem*{maintheorem*}{Main Theorem}
\newtheorem*{thm*}{Theorem}
\newtheorem*{thma*}{Theorem A}
\newtheorem*{thmaa*}{Theorem A'}
\newtheorem*{thmb*}{Theorem B}
\newtheorem*{thmo*}{Theorem 1.1}
\newtheorem*{thmc*}{Theorem C}
\newtheorem*{thmd*}{Theorem D}
\newtheorem*{thmf*}{Theorem 4.1}
\newtheorem*{remark*}{Remark}
\newtheorem*{conjecture*}{Conjecture}
\newtheorem*{prop*}{Proposition}
\newtheorem*{lem*}{Basic Lemma}
\newtheorem{thm}{Theorem}[section]
\newtheorem{cor}[thm]{Corollary}
\newtheorem{lem}[thm]{Lemma}

\theoremstyle{definition}

\newtheorem*{proofc*}{Proof of Theorem C}

\def\vol{\rm{vol}}

\def\SL{\rm{SL}}
\def\rank{\rm{rank}}

\def\hh{\hspace{0.5mm}}

\begin{document}

\author{Youssef Lazar}
\email{ylazar77@gmail.com}
\date{}

\title{On the density of $S$-adic integers near some projective $G$-varieties}

\maketitle

\begin{abstract}  We provide some general conditions which ensure that a system of inequalities involving homogeneous polynomials with coefficients in a $S$-adic field has nontrivial $S$-integral solutions. The proofs are based on the  strong approximation property for Zariski-dense subgroups and  adelic geometry of numbers. We give two examples of applications for systems involving quadratic and linear forms. 
 \end{abstract}


\label{intro}

 \medskip
 \section{Introduction}
 

 Given a finite set of places $S$ of $\mathbb{Q}$ which contains the archimedean one, we consider a finite family of homogeneous polynomials $(f_{i,p})_{p\in S} (1\leqslant i \leqslant r)$ where each $f_{i,p}$ has coefficients in the completion of $\mathbb{Q}$ relative to the place $p\in S$. We are interested in the following problem, given any real $\varepsilon >0$,  can we find an nonzero vector $x\in \mathbb{Z}_{S}^{n}$ such that 
\begin{equation}\label{ineq}0< |f_{i,p}(x)|_{p} \leq \varepsilon \hspace{0.5cm} \mathrm{for}~~\mathrm{every} ~~p\in S  \hspace{0.2cm} \mathrm{and}  \hspace{0.2cm}i=1,\ldots, r ~?
\end{equation}
\noindent Despite its apparent simplicity, this question is extremely difficult to solve in general and as far as we know, only few cases have been settled. Our point of departure is the case of a single isotropic quadratic form $f_1=Q$ for which a solution was found only quite recently by Borel and Prasad \cite{BP} for $S=S_{\infty}$ and completed in the general case as soon Ratner gave a complete solution to the Raghunathan conjecture \cite{bor} in full generality. Their result is an $S$-arithmetic generalization of Margulis' proof of the Oppenheim conjecture \cite{mar}. For the reader interested in such dynamical methods and the applications of Ratner's theory to number theoretical problems we refer to (e.g. \cite{KSS}.).\\
The main tool we are going to use in order to treat the question (\ref{ineq}) with the highest level of generality is the \textit{strong approximation property for algebraic groups}. In other words, we will be merely focusing on  the arithmetical properties of groups actions rather than their ergodic behaviour. It is not very surprizing that strong approximation could solve the same density problems as Ratner's orbit closure theorem does since both results take place in groups generated by one-dimensional unipotent elements. The bridge between these two notions is the Kneser-Tits conjecture which asserts that any simply-connected group which is simple and isotropic over a local field, is generated by its one dimensional unipotent elements which was solved by Platonov for such groups (see e.g. \S 7.2, \cite{PR}). Our main purpose is to provide a substitute in the case when Ratner's theorem fails to hold.  
 One illustration of the great advantage of using strong approximation rather than Ratner's theory is that we get rid of the heavy task of classifying intermediate Lie groups, indeed in most case can such classification is unfeasible unless in the rare cases when we are able to reduce to lower dimension. Fortunately this is the case of the original proof of the Oppenheim conjecture which has been proved for $n=3$ and for Gorodnik's result for pairs $(Q,L)$ which reduces to the dimension four but for pairs the classification of intermediate is much more involved (see \cite{gorod}).  \\ 
In the same circle of ideas, recently Ghosh, Gorodnik and Nevo developped in (\cite{GGN1}, \cite{GGN2}, \cite{GGN3}, \cite{GGN4}) the metric theory of diophantine approximation on homogeneous varieties of semisimple groups in the $S$-arithmetic setting. Among many other deeps results, they proved analogs of Khintchine's and Jarnik's theorems for $S$-adic homogeneous spaces using both ergodic theory and strong approximation for algebraic groups combined with deep concepts coming from the theory of automorphic forms and representation theory. Their method also provides a quantitative version of the strong approximation theorem in homogeneous spaces of semisimple groups.

Finally one should mention that for higher degrees, i.e. when the number of variables and the degrees of the $f_{i}$'s are greater than the number of $r$ of polynomials, the circle method of Hardy and Littlewood still remains the most powerful method for proving (\ref{ineq}) in great generality. In fact, it is providing also sharp quantitative estimate of the numbers of solutions with bounded heights of the number of points lying exactly on a variety following the program promoted by Y. Manin in the seventies.

\subsection{Background and notations}  Let us denote by $\Sigma_{\mathbb{Q}}$ the set of all places in $\mathbb{Q}$, these are given by the set of all prime numbers and the archimedean place corresponding to $\infty$. Let $S$ be a finite set of places in $\Sigma_{\mathbb{Q}_{p}}$ which contains the archimedean one, and let us denote by $S_f$ the subset of all finite (prime) places in $S$, thus we have $S=S_{f} \cup \{\infty\}$. For each prime $p$,  we can define the $p$-adic absolute value is denoted by $|.|_p$  over $\mathbb{Q}$ and we denote by $\mathbb{Q}_{p}$  the corresponding completion of $\mathbb{Q}$. The product of the $\mathbb{Q}_p $ $(p\in S)$ is denoted by $\mathbb{Q}_S$. The set of $p$-adic integers is denoted by $\mathbb{Z}_p$ and is defined to be the set of $x \in \mathbb{Q}$ such that $\vert x \vert_{p} \leq 1$. The ring of $S$-integers of $\mathbb{Q}$ is the set  $\mathbb{Z}_S$ which elements are integral outside $S$ i.e. such that $x\in \mathbb{Z}_p$ for $p \notin S$. For each $p\in S_f$, $\mathbb{Q}_{p}$ is a locally compact (additive) group, hence it is equipped with a Haar measure characterized by the formula $\mu_{p}(a\Omega_p) = |a|_p \mu_p(\Omega_p)$ for all $a \in \mathbb{Q}_p$  and $\Omega_p $ is a measurable subset of $\mathbb{Q}_p$ of finite measure. We normalize it by prescribing the value of the measure $\mu_p$ over the basis of open sets in $\mathbb{Q}_p$ by taking $\mu_{p}(a+p^{n}\mathbb{Z}_p)=p^{-n}$, in particular, $\mu_{p}(\mathbb{Z}_p)=1$. The set of adeles  $\mathbb{A}$ of $\mathbb{Q}$ is the subset of the direct product  $\prod_{p} \mathbb{Q}$ over all the places of $\mathbb{Q}$ consisting of those $x =(x_p)$ such that $ x \in \mathbb{Z}_p$ for almost all places.  The set of adeles $\mathbb{A}$ is a locally compact ring with respect to the adele topology given by the base of open sets of the form $ \prod_{p \in S} \mathbb{Q}_p  \times  \prod_{p\notin S} \mathbb{Z}_p$ where $S \subset \Sigma$ is finite with $ S \supset S_{\infty}$. 
 For any finite subset $S \subset \Sigma$  with $ S \supset S_{\infty}$, the ring of $S$-integral adeles is defined by: 
 \begin{center}
 $\mathbb{A}_{S}  =   \prod_{p \in S} \mathbb{Q}_{p}   \times  \prod_{p \notin S}  \mathbb{Z}_{p}$, thus we can see that $\displaystyle \mathbb{A}=  \bigcup_{S \supset S_{\infty}} \mathbb{A}_S$.
\end{center}

\noindent By definition $\mathbb{Z}_{S} = \mathbb{A}_{S} \cap \mathbb{Q}$, in addition, it can be proved that $\mathbb{Z}_{S}$ is a lattice in $\mathbb{Q}_{S}$ i.e. a discrete subgroup of finite covolume. We can also realize $\mathbb{Z}_{S}$ as a cocompact lattice in $\mathbb{A}_{S}$. One of most fundamental result in arithmetic is the following fact which says that for any nonempty set of places $S$, the image of $\mathbb{Q}$ under the  diagonal embedding  is dense  in $\mathbb{A}_{S}$. This property is called \textit{the strong approximation} for the field $\mathbb{Q}$. For a brief review about this property in the framework of algebraic groups, we invite the reader to read the recent account about this question in (\cite{rapin}), for more details we advice one the most complete reference about this topic (\cite{PR}). \\

\noindent \textit{Quadratic forms over local fields}.~  A quadratic form in $n$ variables over a local field $k$ (i.e of characteristic zero) is given by a symmetric bilinear form $B$ over $k$ such that $Q(x)=B(x,x)$ for any $x\in k^{n}$. We say that $Q$ is \textit{nondegenerate} in $k^n$ if the rank of the matrix associated to $B$ has maximal rank.   A quadratic form $Q_{s}(x)$ is isotropic over $\mathbb{Q}_{s}$ if there exists a nonzero vector $x$ such that $Q_{s}(x)=0$. Over $S$-adic products, a quadratic $Q=(Q_{s})_{s\in S}$ over $\mathbb{Q}_{S}$ is said to be nondegenerate (resp. isotropic) if and only if $Q_{s}$ is nondegenerate (resp. isotropic) over $\mathbb{Q}_{s}$ for each $s\in S$. The special orthogonal group of a quadratic form $(Q_{s})_{s\in S}$ is the product of $S$ of orthogonal groups $\mathrm{SO}(Q_{s})$, the latter is a Lie group which is semisimple as soon as $Q_{s}$ is nondegenerate. The orthogonal group $\mathrm{SO}(Q_{s})$ is isotropic over $\mathbb{Q}_{s}$ (i.e.  it has no nontrivial characters over $\mathbb{Q}_{s}$) if and only if $Q_{s}$ is isotropic. It is well-known that over local fields,  $\mathrm{SO}(Q_{s})$ is isotropic if and if it is noncompact. If $H = \prod_{s\in S} H_{s}$ is a product of $s$-adic Lie groups and $S_{1} \subseteq S$ be a finite subset of places, then $H$ is said to be isotropic over $S_{1}$ if for every $s\in S_{1}$, $H_{s}$ is isotropic over $\mathbb{Q}_{s}$. A quadratic form $(Q_{s})_{s\in S}$ is said to be (globally)\textit{ rational} if there exists a form $Q_{0}$ with rational coefficients such that $Q= \lambda Q_{0}$ for some nonzero $\lambda \in \mathbb{Q}_{S}$, and \textit{irrational} otherwise. A crucial fact is that a quadratic form $Q = (Q_{s})_{s\in S}$ can be irrational (globally) with $Q_{s}$ being rational at some place $s\in S$. Any quadratic space $(V,Q)$ over $k$, can be decomposed in virtue of the Witt's decomposition theorem as follows
$$ V  = V_{an} \oplus \mathrm{rad}(Q) \oplus P_{1} \oplus \ldots \oplus P_r$$  
where the restriction $Q$ to $V_{an}$ is anisotropic over $k$, $\mathrm{rad}(Q)$ is the radical of $Q$ it is equal to zero if $Q$ is nondegenerate and $P_{i} (1\leqslant i \leqslant r)$ are hyperbolic planes such that the restriction of $Q$ to each $P_i$ is isotropic. In particular, $P_{1} \oplus \ldots \oplus P_r$ is the maximal isotropic subspace, $2r$ is called the isotropy index and $r$ is called the Wiit index denoted $i(Q)$, remark that  $Q$ is isotropic over $k$ if and only if $i(Q) \geq 1$.\\

\noindent \textit{Algebraic actions on projective varieties}.~  If we consider an algebraic variety defined over $\mathbb{Q}$ by a prime ideal $I=\langle  f_1,\ldots, f_r \rangle$ for which closed points of given by
$$ X=  \big\{ x \in \overline{\mathbb{Q}}^n  ~|~  f_i (x)=0 , ~~1\leq i \leq r \big\}.$$
The polynomials defining the prime ideal $I$ are supposed to be \textit{homogeneous} in $n$ variables with rational coefficients. For each $1 \leq i \leq r$, if we denote $d_i = \deg f_i$  thus we must have that for any $\lambda \in \overline{\mathbb{Q}}$, $f_{i}(\lambda x)=\lambda^{d_i} f_{i}(x)$, in particular the zero locus  $X=V(I)$  can be seen an \textit{algebraic projective variety}.  Neverthless for practical reasons we need our varieties to be embedded in vectors spaces thefore we will be stuck with the affine point view.
  
\noindent   Let $G$ be the special linear algebraic group $\mathrm{SL}_{n|\mathbb{Q}}$ defined over $\mathbb{Q}$ and let us consider its left action on  the $\mathbb{Q}$-vector space $V= \mathbb{Q}[x_1, \ldots, x_n]$ which is given for each $g \in G$ and $f \in  V$ by 
 \begin{center}
 $ g.f(x) = f(g^{-1}x)$   for all  $x\in \mathbb{Q}^n$.
 \end{center}

\noindent We would like to define a subgroup of $G$  which leaves globally invariant the variety $X=V(I)$ but also every element of the ideal of definition $I$.  For this purpose, we introduce the subgroup $H$ of $G$ defined over $\mathbb{Q}$ given by
 $$H = \big\{ h\in G  ~|~   h.f_{i} =f_{i} , ~~1\leq i \leq r \big\}.$$

\noindent From its very definition  $H$ is an algebraic group which acts linearly on the variety $X$. The\textit{ automorphism group} of $X$ under the action of $G$ ($\mathrm{PSL}_{n}$ in the projective setting) is defined as 
$$\mathrm{Aut}_{G}(X) =\{ x\in X  ~|~ g.x =x\}.$$
It is immediate to see that we have the following inclusions
$$ \mathrm{Aut}_{G}(X) \subseteq H \subseteq G.$$
In particular, one can note that $X$ can have a large group $H$ while having a small or maybe even finite automorphsim group. For this reason,  we prefer to consider the action of $H$ rather than $ \mathrm{Aut}_{G}(X)$ which in the ideal case would be a semisimple noncompact Lie group.

\noindent Let us consider a variety $X/\mathbb{Q}$ embedded diagonally  in a finite product of completions relative to a finite set of places  $S$ in $\mathbb{Q}$ containing the archimedean one. To do so, we first consider the natural diagonal embedding of $\mathbb{Q}$ in $\mathbb{Q}_{S}$ whose image is the direct product of completions $\mathbb{Q}_{p}$ for each $p$ in $S$. Using this embedding we can  consider the family of polynomials $(f_{i,p})_{p\in S}$ over $\mathbb{Q}_{S}$ and for each $p\in S$ we define $X_{p}$ to be the zero set of the $(f_{i,p})$ $(1\leq i \leq r)$, this defines an affine variety over $\mathbb{Q}_{p}$. Therefore $X_{S}$ is the direct product of the completions of $X_p$ ($p\in S$) in $\prod_{s\in S} \overline{\mathbb{Q}}_{p}^{n}$ and we can define in the same way by taking the $S$-product $H_{S}=\prod_{p\in S} H_p$. Clearly the action of $H$ on $X$ induces an equivariant action of $H_S$ on $X_S$ with respect to the diagonal embedding. At some point, we will have to consider the set of $\mathbb{Q}_{p}$-points of $X_p$ which is given by $X_{p}(\mathbb{Q}_{p})= X_p \cap \mathbb{Q}_{p}^{n}$, it is equipped with the $p$-adic topology induced by the base field $\mathbb{Q}_{p}$.
\noindent For each $p\in S$ the set $X_{p}(\mathbb{Q}_{p})$ of  $\mathbb{Q}_{p}$-points of  $X_p$ can endowed with a structure of analytic variety over $\mathbb{Q}_{p}$. 
Using some analogy with analytic complex geometry,  given any real $\varepsilon >0$ we define the $\varepsilon$-\textit{tubular }\textit{neighborhood }of $X_S$ for the $S$-adic topology by 

\begin{center}
$X_{S}^{\varepsilon} = \big\{x\in \mathbb{Q}_{S}^{n} \hh \hh : \hh  | f_{j,s}(x)|_{p}  \hh  \leq \varepsilon~~~~$ for every $ 1 \leq j \leq r $ and every $p\in S \big\}.$
\end{center}

\noindent The elements of $\mathbb{Q}_{S}^n$ which are in $X_{S}^{\varepsilon}$ are called $\varepsilon$-near vectors to $X_S$.  The fact that we have chosen the $f_i$'s to be homogeneous implies that for any  $\varepsilon >0$, the intersection of $ \displaystyle X_{S}^{\varepsilon} $ with the lattice $  {\mathbb{Z}}_{S}^{n} $ contains at least the null vector. Thus for any $\varepsilon >0$, showing that  $\displaystyle X_{S}^{\varepsilon} \cap  {\mathbb{Z}}_{S}^{n} \neq \big\{0\big\}$
amounts to find a nonzero $x\in \mathbb{Z}_{S}^{n}$ such that for every $p\in S$ and $1 \leq i \leq r$ we have, 
\begin{center}
$| f_{i,p}(x) |_p \leq \epsilon.$ 
\end{center}
If we consider the polynomial map $ \phi: \mathbb{Q}_{S}^n $  $ \longrightarrow$   $\mathbb{Q}_{S}^r $ associated to $X_S$, given by

 \begin{center}
$  \phi( x)= \left( f_{1}(x), \ldots, f_{r}(x)\right).$
\end{center}

\noindent then  $\displaystyle X_{S}^{\varepsilon} \cap  {\mathbb{Z}}_{S}^{n} \neq \big\{0\big\}$  if and only if  the origin in $\mathbb{Q}_{S}^r$ is an accumulation point in $\phi(\mathbb{Z}_{S}^{n} )$ for the $S$-adic topology.  A more ambitious question regarding our initial problem  in (\ref{ineq}) is to ask the density of $\phi(\mathbb{Z}_{S}^{n} )$ in $\mathbb{Q}_{S}^r$, note that for real places (\ref{ineq}) , i.e. non-discreteness at the origin, is equivalent to density.\\

\noindent  \textit{Notations.}
\vspace{0.3cm}

\noindent  $\bullet$ If $A$ is a subset of $\mathbb{Q}$, we denote by $\overline{A}^{(p)}$ its $p$-adic closure in $\mathbb{Q}_{p}$ and for any set of places $S$ in $\Sigma_{\mathbb{Q}}$,  ~$\overline{A}^{(S)} = \prod_{p\in S} \overline{A}^{(p)}$. \\

\noindent $\bullet$ If $A$ and $B$ are two sets, we denote $A+B$ by 
 $$A+B = \{ a+b ~|~  a\in A ~\mathrm{and}~~ b \in B\}.$$
\noindent More generally if $A_1, \ldots, A_l$ are sets, their Minkowski sum is denoted 
$$A_1 + \ldots +A_l = \sum_{1\leq i \leq l}^{(M)} A_i  =  \{ a_1+ \ldots + a_n ~|~  a_i\in A_i ~\mathrm{for}~~ 1\leq i \leq l \}.$$
In particular for every  integer $n>0$, the $n$-times sum $A+ \ldots +A$ is denoted by $n\ast A$.

\noindent  $\bullet$ For each $s\in S$, we denote by $\mathrm{Sym}^{2}(\mathbb{Q}_{s}^{n})$ the set of bilinear symmetric forms with coefficients in $\mathbb{Q}_{s}$ with $n$ variables.  We identify this set with the set of all quadratic forms in $n$ variables with coefficients in $\mathbb{Q}_{s}$ since we are in characteristic zero.\\

 \subsection{Main results}The main result gives sufficient conditions in order to ensure that the system (\ref{ineq}) has nontrivial solutions. These conditions are realized if we can find a \textit{rational triple} $(H_0,X_{0}, f_0)$ where $H_0$ acts on $X_{0}=V(f_{0})$ and where $H_{0}$ satisfies some assumptions prior to the application of the strong approximation property. As an application, we provide an answer to (\ref{ineq}) for systems involving one quadratic form and one/several linear form(s) which are well-understood in the real case, i.e.  $S=\{\infty\}$. The three dimensional case was treated by S.G. Dani and G.A. Margulis in \cite{DM1} where solutions to (\ref{ineq}) was provided for pairs $(f_1,f_2) =(Q,L)$. In higher dimensions a similar result has been proved by A. Gorodnik \cite{gorod}. The case of values of quadratic forms restricted to affine subspaces defined linear forms  has been treated by S.G. Dani \cite{D1}. Very recently by the dual case concerning the values of linear forms on a quadric hypersurfaces, O. Sargent was able to prove an Oppenheim type density result in the real case \cite{sargent}. All these results relies on deep results from the ergodic theory of unipotent flows on homogeneous spaces. The most powerful tool to prove such density results is Ratner's orbit closure theorem. The Oppenheim conjecture is a direct consequence of the Raghunathan conjecture, but this conjecture was not yet proved at the time when Margulis and then Borel-Prasad published their results. The Raghunathan conjecture was proved by M.Ratner in full generality and also for $p$-adic Lie groups, the latter has been proved by different methods by Margulis and Tomanov \cite{MT}.
This was fortunate, at least for our work, that Borel and Prasad did not have Ratner's on hand. This has forced them to coin a very interesting method in (\S 4, \cite{BP}) which proves (\ref{ineq}) for a quadratic form with coefficients in $\mathbb{Q}_{S}$ in $n\geq 3$ variables. Under the assumption that $Q$ is \textit{globally irrational} unless at some place  $v\in S$ where $Q_v$ is rational, Borel and Prasad proved  that (\ref{ineq}) holds. This \textit{local rationality} condition is absolutely paramount in order to apply the \textit{strong approximation property} to the orthogonal group of $Q$.  The theorem \ref{main} is generalization of their method in the most general setting so that we can ensures solutions of (\ref{ineq})  for  homogeneous polynomials of arbitrary degrees.  To give some credit to the point of view made by Borel and Prasad, we apply our main theorem to varieties of the form $V(Q,L_0)$ and $V(Q,L_1, \ldots, L_r)$ where $L_0, l_1, \ldots, L_r$ are linear forms. Indeed for both varieties we are able to give sufficient conditions so that (\ref{ineq}) holds. We expect that it holds for more general examples of $G$-varieties.\\

  \begin{thm}\label{main} Let $S$ be a finite set of places in $\mathbb{Q}$ containing the archimedean one. For each $s\in S$, we are given a projective algebraic variety $X_{s}$ over $\mathbb{Q}_s$ defined by a homogeneous prime ideal $I_s$ of $\mathbb{Q}_{s}[x_1, \ldots, x_n]$ and let $H_s$ be the algebraic $\mathbb{Q}_s$-subgroup of $\SL_{n}(\mathbb{Q}_s)$ leaving invariant every generator of $I_s$. Assume that the following subset of places  $S_{1} \subset S_{f}$   is nonempty, $$S_{1} = \{ p \in S_{f} ~|~   H_{p} ~~\mathrm{is} ~~\mathrm{rational}~~\mathrm{over}~~\mathbb{Q} \}.$$
If there exists a connected algebraic subgroup $H_{0}$ of $G$ rational over $\mathbb{Q}$ and a hypersurface $X_{0}=V(f_0)$ defined over $\mathbb{Q}$ such that
\begin{enumerate}
\item $H_0$ is a semisimple absolutely almost simple algebraic  $\mathbb{Q}$-group.

\item For every prime $p\in S_1$, $X_{p}=X_0$ and $H_{p}=H_{0}$.

\item  Every $\mathbb{Q}$-simple factor of $H_0$ is isotropic over $S_{1}$.

\end{enumerate} 
 
 \noindent Then for any $\epsilon >0$, there exists a nonzero $S$-integral vector lying $\epsilon$-near $X_S$ i.e. 
$$\displaystyle X_{S}^{\varepsilon} \cap  {\mathbb{Z}}_{S}^{n} \neq \big\{0\big\}.$$

 \end{thm}
 
As an application we can use the theorem  to prove existence of near integral vectors to the variety $X_S=\{Q=L=0\}$  which can be see as the hypersurface consisting of one nondegenerate quadric $Q=0$ cutted out by the hyperplane of  equation  $\{ L=0\}$. In the following result we extend the validity of a previous result of the author(Corollary 2.2, \cite{yl}) proved under the condition that any nontrivial linear combinaison $ \alpha_s Q_s+ \beta_s  L_{s}^{2}$ should be irrational at all places $s\in S$. Indeed we are able to prove that the same result holds when if we only assume that  $ \alpha Q+ \beta  L^{2}$  is (globally) irrational over $\mathbb{Q}_{S}$ leaving the possibility that it is rational at some  nonarchimedean place. 
\begin{cor}\label{cor-pairs} Assume $S$  is a finite set of places in $\mathbb{Q}$ containing the archimedean one and let $Q = (Q_s)_{s \in S}$  be a quadratic form  and $L = (L_{s})_{s \in S}$ be a linear form on $\mathbb{Q}_{S}^{n}$ with $n \geq 4$ and $L_s \neq 0$ for all $s\in S$. Suppose that the pair $(Q,L)$ satisfies the following conditions, 
 
  \begin{enumerate}
\item $Q$ is nondegenerate.
\item $Q_{\vert L = 0}$ is nondegenerate and isotropic.
\item For any choice $\alpha,\beta$ in $\mathbb{Q}_{S}$ with $( \alpha,\beta) \neq (0,0)$,   the form $ \alpha Q+ \beta  L^{2}$ is  irrational with at least  one place $v\in S_f$ where $ \alpha_v Q+ \beta_v  L^{2}$ is proportional to a rational form.
\end{enumerate}
Then for any $ \varepsilon>0$,  there exists a nonzero $x \in {\mathbb{Z}}_{S}^{n}$ such that  
  \begin{center}
 $ |Q_{s}(x)|_s <  \varepsilon $  and $ |L_{s}(x)|_{s} <  \varepsilon $ ~~ for each $s \in S$.
 \end{center} 
\end{cor}
\vspace{0.5cm}
\noindent In the same vein, we are also able to treat the case when we  consider several linear forms instead of one. This result is a variation of a recent due to O. Sargent \cite{sargent} in the $S$-arithmetic setting which is proved using Theorem \ref{main}.

\begin{cor}\label{cor-pairs-2} Assume $S$  is a finite set of places in $\mathbb{Q}$ containing the archimedean one and suppose that $Q$ is an isotropic nondegenerate quadratic form in $n\geq 3$ variables over $\mathbb{Q}_{S}$ with rational coefficients. Let  $M=(L_1, \ldots, L_{r})$ be a linear map $M: \mathbb{Q}_{S}^{n} \rightarrow \mathbb{Q}_{S}^{r}$ where $(L_{i,s})_{s \in S}$ $(1\leqslant i \leqslant r)$  are linear forms of rank $r$ over $\mathbb{Q}_{S}^{n}$ in $n \geq 3$ variables which satisfies the following conditions:
 
  \begin{enumerate}
\item $n> \max_{s\in S} (\dim \ker M_{s} ) + 2$.
\item  \emph{rank} $(Q_{s \vert  \mathrm{Ker} M_s}) =r$ and $Q_{ \vert \mathrm{Ker} M_s}$ is isotropic and , ~for every $s\in S$. .
\item For each choice of $\alpha_1,\ldots, \alpha_r$ in $\mathbb{Q}_{S}$ with $( \alpha_1,\ldots, \alpha_r) \neq (0,\ldots, 0)$,   the linear form $ \alpha_{1} L_{1}+ \ldots +  \alpha_{r} L_{r}$ is  irrational unless at one place $v\in S_f$ where it is proportional to a rational form.
\end{enumerate}
Then for any $ \varepsilon>0$,  there exists a nonzero $x \in \mathbb{Z}^{n}_{S}$ such that  
  \begin{center}
$Q_{s}(x)=0$ and $ |L_{i,s}(x)|_{s} <  \varepsilon $ ~~ for each $s \in S$ and $1\leqslant i \leqslant r$.
 \end{center}

\end{cor}

\vspace{0.5cm}

\subsection*{Remarks} 
$(1)$  The proof of Theorem \ref{main} is based on a strenghtening of the strong approximation theorem which apply to Zariski dense subgroups of reductive groups proved by Matthews, Vaserstein and Weisfeiller (\cite{matt}, \cite{weis}) and later by Nori \cite{nori}. In the meanwhile, Venkataramana (Proposition $(5.3)$, \cite{venky}) proved that there exists Zariski dense subgroup of integral points of $SL_{n}(\mathbb{Z})$ which contains no unipotent elements and thus such subgroup might be eligible for strong approximation even in it is far from being unipotent or generated by unipotents elements. For more general details about strong approximation in algebraic groups and more particularly this version, we advice the reader the recent survey of A.S. Rapinchuk \cite{rapin}.\\

\noindent $(2)$ The method used in Theorem \ref{main} is an adaptation of the work of Borel-Prasad in the case when $S$ contains at least one place where the form is rational (see \S 4, \cite{BP}). \\

 \noindent $(3)$ A nice feature of the proof of Corollaries \ref{cor-pairs} and \ref{cor-pairs-2} is that we do not have to reduce to lower dimension. Indeed, the reduction process was a pre-condition for proceeding to the classification of intermediate  subgroups aring from the application Ratner's theorem. \\
 
  \noindent $(4)$ In the statement of Corollary \ref{cor-pairs-2} we have not tried to find optimal and thus it can be seen as a partial $S$-adic generalization of the main result of O. Sargent in the real case \cite{sargent}. It may be possible to refine the conditions but our aim was to give an easy example of application of the main theorem.  \\




\section{Proof of Theorem \ref{main}}
We denote by $S_2$ the set of places of $S$ which are disjoint from $S_1$, in particular $S_1$ contains only nonarchimedean places  with let us say $S_1 = \{ p_1, \ldots, p_s \}$ and $S_{2}=\{q_1, \ldots, q_{l-1}, \infty \}$.
Let us define $\Lambda$ to be the stabilizer of the standard lattice 
 $\mathbb{Z}_{S}^{n}$ in $H_{0}(\mathbb{Q})$ i.e. 
 $$\Lambda = \{ h \in H_{0}(\mathbb{Q}) ~|~   h(\mathbb{Z}_{S}^{n} ) = \mathbb{Z}_{S}^{n} \}.$$
 Under the diagonal embedding, $\Lambda$ can be seen an $S$-arithmetic which is discrete in $H_{0}(\mathbb{Q}_{S})$. Consider the universal $\mathbb{Q}$-isogeny $\pi: \tilde{H_{0}} \rightarrow H_{0}$, here $\widetilde{H_{0}}$ is a semisimple simply connected group defined over $\mathbb{Q}$. Let us choose an arbitrary $S$-arithmetic subgroup $\widetilde{\Lambda}$ of $\tilde{H_0}(\mathbb{Q})$  embedded  as a discrete subgroup in $\tilde{H_{0}}(\mathbb{Q}_S)$. Since $\tilde{H_0}$ is absolutely almost simple over $\mathbb{Q}$ and every simple factor $H_i$ of $\widetilde{H_{0}}$ is noncompact at every place in $S_1$, in particular $\sum_{s\in S} {\rank_{\mathbb{Q}_v}} H_i >0$. Borel's density theorem (see e.g.  proposition 3.2.10 (\cite{Mbook}) applied  to $\tilde{H_0}$ with $K=k=\mathbb{Q}$ gives us  that the discrete subgroup $\tilde{\Lambda}$ is Zariski-dense in $\tilde{H_0}$. A special instance of the strong approximation for Zariski-dense subgroups (see e.g. \cite{PR}, Theorem 7.14) applied in $\tilde{H_0}$ for the set of primes $S_{1}=\{p_1, \ldots, p_l \}$ implies that the closure $\overline{\tilde{\Lambda}}$ is open in $\tilde{H_0}(\mathbb{A}_{S_1})$. The latter adelic group is described by the following product $$   \prod_{p \in S_{1} }\widetilde{H}_0(\mathbb{Q}_{p}) \times  \prod_{q \in S_{2}}\widetilde{H}_0(\mathbb{Z}_{q})\times  \prod_{q \notin S_{1} \cup S_{2}}\widetilde{H}_0(\mathbb{Z}_{q})  $$
Thus for every $q \in S_2$, the projection $\overline{\widetilde{\Lambda}}^{(q)}$ is open in $\tilde{H_0}(\mathbb{Z}_{q})$. The universal $\mathbb{Q}$-isogeny $\pi$ transforms $\tilde{\Lambda}$ into an arithmetic subgroup $\pi(\tilde{\Lambda})$ in $H_{0}(\mathbb{Q})$ such that $\overline{\pi(\widetilde{\Lambda})}^{(S_{2})}$ is open in 
$H_{0}(\mathbb{Q}_{S_{2}})$.  The two arithmetic subgroups $\pi(\tilde{\Lambda})$ and $\Lambda$ are commensurable in $H_{0}$, thus $\pi(\tilde{\Lambda}) \cap \Lambda$ is a finite index subgroup in $\pi(\tilde{\Lambda}) $ (see e.g. cor. 3.2.9 \cite{Mbook}). For each $p\in S$, we define $U_{p}$ to be the projection of $\pi(\tilde{\Lambda}) \cap \Lambda$ onto its $p$-component, that is, 
\begin{equation}\label{Up}
U_{p} := (\pi(\tilde{\Lambda}) \cap \Lambda)_{p}.
\end{equation}
As we have seen above, the $p$-adic closure  $\overline{U_p}^{(p)}$ lies in the open subgroup  $\overline{\pi(\tilde{\Lambda})}^{(p)}$  of  $H_{0}(\mathbb{Q}_{p})$ for each $p\in S_{2}$ and $U_{S_{2}}$ is contained in $\overline{\Lambda}^{(S_{2})}$ by definition.\\
Now let us introduce the following subset of vectors in $\mathbb{Q}_{p}^{n}$ defined for each $p\in S_2$, by:
$$ \mathfrak{X}_{p} = U_{p}^{-1} X_{p}(\mathbb{Q}_{p})-X_{0}(\mathbb{Q}_{p})$$
where $U_{p}^{-1}= \{ a_{p}^{-1} | a_{p}\in U_{p} \}$ is an open subset of invertible matrices in $H_{0}(\mathbb{Q}_{p})$ where $U_p$ was defined in (\ref{Up}). \\
We claim that $ \mathfrak{X}_{p}$ is a nonempty open blunt cone in $\mathbb{Q}_{p}^{n}$ for $p\in S_2$. Indeed, since $X_{p}$ is not rational over $\mathbb{Q}$ for each $p\in S_2$ while $X_{0}$ is,  this forces $X_{p}\neq X_{0}$ for $p\in S_2$. The Zariski density of $X_{p}(\mathbb{Q}_{p})$ (resp. $X_{0}(\mathbb{Q}_{p})$)  in $X_{p}$ (resp. $X_{0}$) implies that 
$X_{p}(\mathbb{Q}_{p}) \neq X_{0}(\mathbb{Q}_{p})$ for each $p\in S_2$. In particular for each $p\in S_2$, there exists an $x\in X_{p}(\mathbb{Q}_{p}) - X_{0}(\mathbb{Q}_{p})$, thus taking $a_{p}=I_{n}$ we get that $x \in  \mathfrak{X}_{p}$. Now let us fix $x \in  \mathfrak{X}_{p}$ and consider $y\in \mathbb{Q}_{p}^{n}$ sufficiently close to $x$, so that we can find an $g\in {\SL}_{n}(\mathbb{Q}_p)$ such that $y=gx$ where $g$ is close to $I_n$ in ${\SL}_{n}(\mathbb{Q}_p)$. Since $I_n \in U_p$  then we can assume that $g$ is arbitrarily close to $I_n$ in the open set $U_{p}$.  We have $f_{i,p}(g_{p} x)=0$ for some $g_{p} \in U_{p}$ close to $I_n$ and since $y=gx$ we get $f_{i,p}((g_{p}g^{-1}) y )= f_{i,p}(g_{p}x ) = 0$ with $g_{p}g^{-1} \in U_{p}$ thus $y\in  \mathfrak{X}_{p}$. It is clear that $ \mathfrak{X}_{p}$ does not contains the origin and that for any $\lambda \in \mathbb{Q}_{p}^{\ast}$ and $x\in  \mathfrak{X}_{p}$, $\lambda x \in  \mathfrak{X}_{p}$ by homogeneity of the $f_i$'s, the claim is proved.\\
For each $p\in S_1$, denote by  $\Delta_{p}(r)$ the hypercube centered at $0$ with radius $r$ in $\mathbb{Q}_{p}^{n}$, i.e.
$$ \Delta_{p}(r) =\big\{ x \in \mathbb{Q}_{S_{1}}^{n} ~|~  |x_i|_{p} \leq r ~~  \big\} $$
and denote $\Delta_{S_{1}}(r) = \prod_{p\in S_{1}} \Delta_{p}(r)$. Our aim is to show that $\cap_{\varepsilon >0} X_{S}^{\varepsilon} \cap \mathbb{Z}_{S}^{n} \neq \{ 0\}$, a first step consists to show the existence of nonzero lattice point for the domain $(\Delta_{S_{1}}(\delta) \times \mathfrak{X}_{S_2})$ which can be seen as a sort of approximation of $X_{S}^{\varepsilon}$.

\begin{lem}\label{lemma} For every $\delta > 0$, we have
$$(\Delta_{S_{1}}(\delta) \times \mathfrak{X}_{S_2}) \cap \mathbb{Z}_{S}^{n} \neq \{ 0\}.$$
\end{lem}
 
\noindent \textbf{Proof of the lemma.} Let us fix a place $p\in S_2$  whether it is archimedean or not, and set $l=|S_2|$. Let $v_{1,p}\in \mathbb{Q}_{p}^{n}$ be a nonzero in vector $\mathfrak{X}_{p}$ and complete it to a basis $\{v_{1,p}, v_{2,p}, \ldots, v_{p,n}\}$ of $\mathbb{Q}_{p}^{n}$.  For each real $a>0$ we introduce the hypercubes $V_p(a)$ in $\mathbb{Q}_{p}^{n}$ and  $W_p(a)$  in $\mathbb{Q}_{p}^{n-1}$ defined as
\begin{center}
$V_{p}(a) = \displaystyle \big\{ \sum_{i=1}^{n}  \alpha_{i}  v_{i,p}  ~|~  |\alpha_i |_{p} \leq a \big\}$  and $W_{p}(a) = \displaystyle \big\{ \sum_{i=2}^{n}  \alpha_{i}  v_{i,p}  ~|~  |\alpha_i |_{p} \leq a \big\}$.
\end{center}

\noindent Since we know that  $\mathfrak{X}_{p}$ is open, we can find an infinitesimal  hypercube $v_{1,p} \oplus W_{p}(\alpha)$ for some small enough real $\alpha >0$  so that it is contained in $\mathfrak{X}_p$. However we have the following fact: the resulting infinitesimal hypercube $v_{1,p} \oplus W_{p}(\alpha)$ remains in the cone $\mathfrak{X}_{p}$ if we perform  a translation in the direction of $v_{1,p}$ away from $\mathbb{Z}_{p}$ \footnote{Note that if $p$ is nonarchimedean $\mathbb{Q}_{p}$ is not an ordered field neither even partially, so one has to be careful with the meaning of this assertion. The real meaning is arithmetical rather than geometrical as it can be seen just below.}. In other words, we can find an real positive $\alpha$ small enough so that  for each $p\in S_{2}$ and for any given arbitrary $\eta \in \mathbb{Q}_{p}$ with $|\eta |_{p} >1$ (i.e. $\eta \notin \mathbb{Z}_p$) we have simultaneously 
\begin{equation}\label{inclusionaway}
\eta v_{1,p} \oplus W_{p}(\alpha) \subset \mathfrak{X}_p ~~\mathrm{and}~~ V_{q}(\alpha)\subset \mathfrak{X}_p.
\end{equation}
Indeed, let us set $u = \eta v_{1,p} +  \sum_{i=2}^{n} \alpha_i v_{i,p} \in \eta v_{1,p} \oplus W_{p}(\alpha)$ for $\eta \in \mathbb{Q}_{p}^{\ast}-\mathbb{Z}_{p}$. Thus $u$ can be written as
$$ u = \eta(v_{1,p}+ \sum_{i=2}^{n} \eta^{-1}\alpha_i v_{i,p} ).$$
It is clear that for every $2 \leq i\leq n$, $|\eta^{-1}\alpha_{i}|_p < \alpha$, thus 
$$v_{1,p}+ \sum_{i=2}^{n} \eta^{-1}\alpha_i v_{i,p} \in v_{1,p} \oplus W_{p}(\alpha) $$
 Therefore using the cone invariance for $ \mathfrak{X}_p$,   we infer that $u \in \mathfrak{X}_p$, which proves the claim (\ref{inclusionaway}).
For this choice of $\alpha$ and for each reals $\delta,t >0$ let us introduce the $S$-adelic domain 

$$\mathcal{C}_{p}(\delta, t) = \Delta_{S_1}(\delta/2l) \times [0,t]v_{1,p} \oplus W_p(\alpha/2) \times  \prod_{q\in S_{2}\backslash\{p\}}   V_q(\alpha/2l) \times \prod_{s\notin S} \mathbb{Z}_{p} \subset \mathbb{A}_{S}^{n}.$$

\noindent Let us define $C_{p}(\delta)$ to be $2\ast \mathcal{C}_{p}(\delta, 1)$. It is a compact subset of $\mathbb{A}_{S}^{n}$ and thus it meets the discrete subgroup $\mathbb{Z}_{S}^{n}$ in finitely many points at the number of $k=|C_{p}(\delta) \cap \mathbb{Z}_{S}^{n}|$. The set of $S$-integral vectors $\mathbb{Z}_{S}^n$ is a cocompact lattice in $\mathbb{A}_{S}^n$ i.e. $\mathbb{A}_{S}^{n}/\mathbb{Z}_{S}^{n}$ is a compact space of finite volume for the measure $\mu$ induced by $\vol_{\mathbb{A}_S}$ on the quotient space. Let us denote by $\mathcal{D}$ a fundamental domain for the quotient $\mathbb{A}_{S}^{n}/\mathbb{Z}_{S}^{n}$ and by $\mu(\mathcal{D})$ its volume, in particular we have 
$$\mathbb{A}_{S}^{n} = \bigcup_{y\in \mathcal{D}} y + \mathbb{Z}_{S}^{n}. $$

\noindent We can find some $\tau$ large enough so that exists $k+2$ disctinct points $y_{0}, \ldots, y_{k+1}$ in $\mathcal{C}_{p}(\delta, \tau/2)$ such that $y_{i} -y_{j} \in \mathbb{Z}_{S}^{n}$ for any $1 \leqslant i<j \leqslant k+1$. Indeed, for this it suffices to remark that the function $t \mapsto {\vol}_{\mathbb{A}_{S}}({\mathcal{C}}_{p}(\delta, t))$ in increasing, so for  $\tau$ large enough we have the inequality 
$$  {\vol}_{\mathbb{A}_{S}}({\mathcal{C}}_{p}(\delta, \tau/2)) > (k+1) \mu(\mathcal{D}).$$
Now applying the adelic Blichfeldt's principle in $\mathbb{A}_{S}^{n}$ (see e.g. Lemma 4, \S 5.2, \cite{m1}) we obtain the required $y_{i}$'s. Let us put $K=\mathbb{A}_{S}^{n}- \mathcal{C}_{p}(\delta)$ and $x_i = y_{0}-y_{i}$ for $1\leq \leq k+1$. The previous assertion tells us that $x_1, \ldots, x_{k+1}$ are nonzero elements of $\mathbb{Z}_{S}^{n}$, namely we have a set of $k+1$ integral vectors. But there are only $k$ elements in $ \mathcal{C}_{p}(\delta) \cap \mathbb{Z}_{S}^{n}$, thus exactly one of them must lie in $K_{p} \cap \mathbb{Z}_{S}^{n}$, call it $x(p)$. By definition $x(p)$ is the difference of two elements in $ \mathcal{C}_{p}(\delta, \tau/2)$, thus $x(p) \in 2 \ast \mathcal{C}_{p}(\delta, \tau/2)$.

\noindent Hence for each $p\in S_2$, we have obtained an nonzero element $x(p)$ of $\mathbb{Z}_{S}^{n}$ which lies in $ K_p \cap   2\ast \mathcal{C}_{p}(\delta, \tau/2)$, the latter subset satisfying the following inclusion 
$$ K_p \cap   2\ast \mathcal{C}_{p}(\delta, \tau/2)\subset \Delta_{S_1}(\delta/l) \times (1,\tau]v_{1,p} \oplus W_p(\alpha) \times  \prod_{q\in S_{2}\backslash\{p\}}   V_q(\alpha/l) \times \prod_{s\notin S} 2\ast \mathbb{Z}_{p}.$$
\noindent  Finally let us set  $x= \sum_{ p \in S_2} \pi_{S}(x(p))$ where $\pi_S$ is the projection to the $S$-factor. It is already clear that $x$ is a nonvector vector in $\mathbb{Z}_{S}^{n}$. It remains to verify that $x \in  \Delta_{S_{1}}(\delta) \times \mathfrak{X}_{S_2}.$
For the $S_1$-components, we have 
$$\pi_{S_{1}}(x)= \sum_{p\in S_{2}} \pi_{S_{1}}(x(p)) \in  \sum_{p\in S_{2}}^{(M)} \pi_{S_{1}}\left( K_p \cap   2\ast \mathcal{C}_{p}(\delta, \tau/2) \right) \subset l \ast \Delta_{S_{1}}(\delta/l) \subset  \Delta_{S_{1}}(\delta).$$
On $S_2$ side, we  isolate the diagonal component in order to obtain
$$\pi_{S_2}(x) = x(p)_{p} +  \sum_{q\in S_{2}\backslash\{p\}} x(p)_{q} \in (1,\tau]v_{1,p} \oplus W_p(\alpha)  \times \prod_{q\in S_{2}\backslash\{p\}}  \sum_{p^{\prime}\in S_{2}\backslash\{p\}}^{(M)}   V_q(\alpha/l)_{p^{\prime}}.$$
Remembering the choice of  $\alpha>0$  made above in (\ref{inclusionaway}), we infer that
$$ \pi_{S_2}(x) \in \mathfrak{X}_p \times  \prod_{q\in S_{2}\backslash\{p\}} V_q(\alpha) \subset \mathfrak{X}_{S_{2}}.$$
Hence for any $\delta>0$, we can always find a nonzero vector $x$ lying in $ (\Delta_{S_{1}}(\delta) \times \mathfrak{X}_{S_2}) \cap \mathbb{Z}_{S}^{n}$ and this achieves the proof of the lemma.
\begin{flushright}
$\square$
\end{flushright}

\begin{figure}
\fbox{\begin{minipage}{15cm}
\begin{tikzpicture}[scale=1.5]
\draw (0.5,1.25) -- (1.5,1.25) -- (1.5,0.4) --(0.5,0.4)--cycle;

\draw[thick,densely dotted, blue,rounded corners=1cm] (0,0)--(1,0.5)--(5,-0.5)--(6,1);
\draw[thick,densely dotted, rounded corners=1cm] (0,0.5)--(1,1)--(5,0)--(6,1.5);
\draw[thick,densely dotted, red ,rounded corners=1cm] (0,0.75)--(1,1.25)--(5,0.25)--(6,1.75);
\draw[thick,densely dotted, blue,rounded corners=1cm] (0,1)--(1,1.5)--(5,0.5)--(6,2);

\node[label=below:$\Delta_p(\delta)$] (A) at (1,0.4) {};
\node[label=below:$0$] (B) at (1,0.9) {};
\node[label=below:$\color{red}{\gamma_p x_p}$] (C) at (3,0.56) {};
\node[label=below:$\color{red}{ x_p}$] (C) at (1,1.9) {};
\node[label=right:$\color{blue}{X_{0}^{\varepsilon}}$] (A) at (6,1) {};
\node[label=right:${X_p =X_0= \{ f_0= 0\}}$] (A) at (6,1.5) {};

\node[label=right:$\Delta_{S_{1}}(\delta) \subset \mathbb{Q}_{S_{1}}^n$] (A) at (8,0) {};
\fill [red] (1,1.1) circle (0.05);

\fill [red] (3,0.75) circle (0.05);

 \fill [black] (1,0.85) circle (0.05);
\end{tikzpicture}

\end{minipage}}
\vspace{1cm}
\fbox{\begin{minipage}{15cm}
\begin{tikzpicture}[scale=1.5]
\draw (2.5,0.9) -- (3.5,0.9) -- (3.5,0.15) --(2.5,0.15)--cycle;

\draw[thick,densely dotted, blue,rounded corners=1cm] (0,0)--(1,0.5)--(5,-0.5)--(6,1);
\draw[thick,densely dotted, rounded corners=1cm] (0,0.5)--(1,1)--(5,0)--(6,1.5);
\draw[thick,densely dotted, blue,rounded corners=1cm] (0,1)--(1,1.5)--(5,0.5)--(6,2);
\fill [black] (3,0.5) circle (0.05);
\fill [red] (3,0.75) circle (0.04);
\fill [black] (5,2.5) circle (0.05);
\fill [black] (1,0.9) circle (0.05);
\node[label=below:$0$] (B) at (1,0.9) {};
\node[label=above:] (G) at (3,0.5) {};
\node[label=above:$u_p x_p$] (T) at (3,-0.05) {};
\node[label=above:$\color{red}{\gamma_p  x_p}$] (R) at (3,0.75) {};
\node[label=left:${x_p}$] (X) at (5,2.5) {};

\draw[->,>=latex] (X) to[bend left=30] (G);

\node[label=right:${X_p = \bigcap_{j} \{ f_{j,p}= 0\}}$] (A) at (6,1.5) {};
\node[label=right:$\color{blue}{X_{S_{2}}^{\varepsilon}}$] (A) at (6,1) {};
\node[label=right:$\mathfrak{X}_{S_{2}} \subset  \mathbb{Q}_{S_{2}}^n$] (A) at (8,0) {};

\end{tikzpicture} 
\end{minipage}}
{\caption{The space between the blue dotted lines is the tubular neighborhood $X_{S}^{\varepsilon}$. }}
\end{figure}
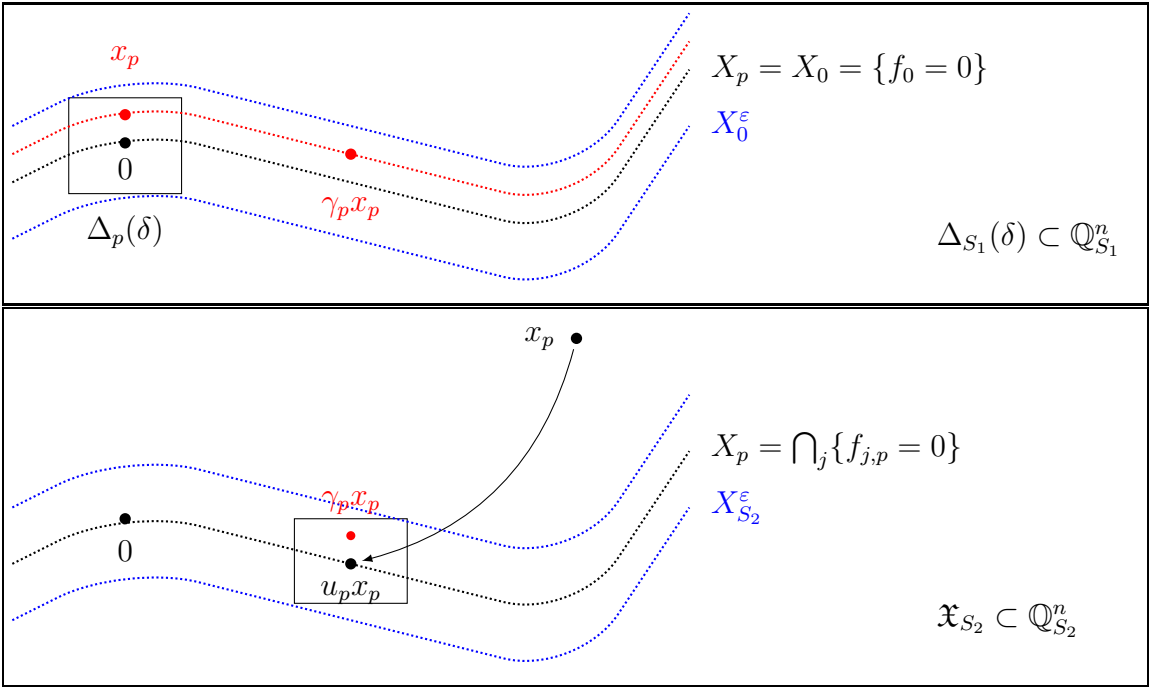

\noindent We are now ready to prove the theorem, for this let us fix an $\varepsilon >0$.  The fact that the $f_i$'s are homogenous polynomials and thus, in particular, continous with $f_i(0)=0$,  implies that the existence of  $\delta(\varepsilon)>0$ such that $\Delta_{S_{1}}(\delta(\varepsilon)) \subset X_{S_{1}}^{\varepsilon}$. Using Lemma \ref{lemma} with $\delta=\delta(\varepsilon)$,  one obtains a nonzero $S$-integral vector $x\in \mathbb{Z}_{S}^{n}$ such that $\pi_{S_{1}}(x) \in X_{S_{1}}^{\varepsilon}$ and $\pi_{S_{2}}(x) \in \mathfrak{X}_{S_{2}}$. The latter condition means that for each $p\in S_{2}$ and corresponding $x_{p}\in \mathfrak{X}_{p}$, there exists some $u_{p} \in U_{p}$ such that $f_{i,p}(u_{p}x_{p})=0.$ Since $U_p$  is open there exists $g_{p} \in U_{p}$ such that $$0 < |f_{i,p}(g_{p}x_{p})|_{p} \leq \varepsilon/2.$$
With the help of the strong approximation theorem, we have seen earlier that $U_{S_{2}}$ is contained in $\overline{\Lambda}^{(S_{2})}$, and  incidentally we find $(\gamma_{p})_{p\in S_{2}} \in \Lambda_{S_{2}}$ such that for every $p\in S_{2}$ and $1\leq i \leq r$, one has
\begin{equation}\label{forS2}
0 < |f_{i,p}(\gamma_{p}x_{p})|_{p} \leq \varepsilon.
\end{equation}

\noindent Consider the projection $\Lambda_{S} \twoheadrightarrow \Lambda_{S_{2}}$ and let $(\tilde{\gamma})_{p\in S}$ be a lift of $(\gamma_{p})_{p\in S_2}$, that is, $\tilde{\gamma}_{p} =\gamma_{p}$ for every $p\in S_{2}$. We claim that $y=\tilde{\gamma}x$ is the solution of our problem. Indeed, on the one hand the inequalities in  (\ref{forS2}) imply that $\pi_{S_{2}}(y) \in X_{S_{2}}^{\varepsilon}$.\\
On the other hand since $\pi_{S_{1}}(x) \in X_{S_{1}}^{\varepsilon}$ and $X_{p}^{\varepsilon}$ is $H_{0}$-invariant with $\Lambda_{p} \subseteq H_0$ for every $p\in S_1$, we deduce that $y_p=(\tilde{\gamma}_{p} x_{p}) \in X_{p}^{\varepsilon}$ for every $p\in S_{1}$, i.e. $\pi_{S_{1}}(y) \in X_{S_{1}}^{\varepsilon}$. This allows us to conclude the existence of a nonzero vector $y \in X_{S}^{\varepsilon} \cap \mathbb{Z}_{S}^{n}$ and this finishes the proof of the  theorem. 

\section{Proof of Corollary \ref{cor-pairs}} Let us consider a pair  $(Q_{s},L_{s})_{s\in S}$ over  $\mathbb{Q}_{S}$ satisfiying all the assumptions of the corollary  \ref{cor-pairs}. For each $s\in S$, we set $X_{s}$ to be the algebraic (projective) variety given by $\{ Q_{s}=L_{s}=0\}$, geometrically this can be seen as the cone $Q_s =0$ cutted out by the hyperplane of equation $L_{s}=0$ in $\mathbb{Q}_{s}^{n}$. It is more suitable here to think $X_s$ as the quadric of equation $ \{ Q_{s\mid_{L_s=0} }=0\}$, in that way the assumptions $(1)$ and $(2)$ amounts to say that this quadric $ \{ Q_{s\mid_{L_s=0} }=0\}$ is nondegenerate and contains at least one nonzero vector in $\mathbb{Q}_{s}^{n}$ for every $s\in S$. We need to introduce the following map associated to the pair $(Q,L)$ over $\mathbb{Q}_{S}$
$$\begin{array}{cccc}
\psi :& \mathbb{P}^{1}(\mathbb{Q}_{S}) & \rightarrow &\mathrm{Sym}^{2}(\mathbb{Q}_{S}^{n})\\
&( \alpha : \beta ) & \mapsto & \alpha Q+ \beta L^2. 
\end{array}$$
This induces at each place $s\in S$, a (local) map $\psi_s(\alpha_s : \beta_s)=  \alpha_s Q_s+ \beta_s L_{s}^{2}$.

\noindent The assumption (3) says that the range of $\psi$ is in the subspace $\mathrm{Sym}_{ir.}^{2}(\mathbb{Q}_{S}^{n})$ consisting
of quadratic forms in $ \mathrm{Sym}^{2}(\mathbb{Q}_{S}^{n}) $ which are not proportional to rational form over $\mathbb{Q}_{S}$.  The complement consisting of quadratic forms over $\mathbb{Q}_{S}^{n}$  (resp. $\mathbb{Q}_{s}^{n}$)  which are proportional to a rational form is denoted $ \mathrm{Sym}_{rat.}^{2}(\mathbb{Q}_{S}^{n})$ (resp. $ \mathrm{Sym}_{rat.}^{2}(\mathbb{Q}_{s}^{n})$).\\

\noindent \textit{Case 1.} If we are in the case where  $\psi_s$ assumes its values in $\mathrm{Sym}_{ir.}^{2}(\mathbb{Q}_{s}^{n})$ for every $s\in S$, there is nothing to prove.  Indeed, the corollary 2.2. in \cite{yl} already gives the required result.\\

\noindent  Now we treat the case of interest here.

\noindent \textit{Case 2.} Suppose there exists a place $v\in S_{f}$ such that the range of $\psi_v$ is in $\mathrm{Sym}_{rat.}^{2}(\mathbb{Q}_{v}^{n})$. This means concretely that for any not both zero couple of constants $(\alpha_v,\beta_v) \in \mathbb{Q}_{v}^2$ we have that $\alpha_v Q_v+ \beta_v L_{v}^{2}$ is proportional to a quadratic form $Q_0$ with rational coefficients, thus  $\alpha_v Q_v+ \beta_v L_{v}^{2} =  \lambda_v Q_0$ for some $\lambda_{v} \in \mathbb{Q}_{v}^{\ast}$, in particular $ \psi_v(1:0)=Q_v$ and $ \psi_v(0:1)=L_v$ are proportional to a rational form. Thus  we  can write $L_v = \mu_v L_0$ where $L_0$ is a rational linear form and $\mu_{v} \in \mathbb{Q}_{v}^{\ast}$  and from the definition of $Q_0$ we have that $\alpha_v Q_{v \mid L_v =0}  =  \lambda_v Q_{0 \mid  L_0 =0 }$ and in particular $SO( Q_{v \mid L_v =0}) = SO( Q_{0 \mid  L_0 =0 })$. Now for each $s\in S$ we denote by $H_s$ the stabilizer of the pair $(Q_s,L_s)$ to be the subgroup of ${\SL}_{n}(\mathbb{Q}_s)$  leaving both invariant $Q_s$ and $L_s$. Following Lemma 4.1 in \cite{yl}, since $Q_v$ and $L_v$  are both  proportional to rational forms  then  we  can arrange a basis $\{w_1, \ldots, w_{n-1},u\}$ consisting of rational vectors with the  condition that  $\{L_v =  0\}=\langle w_1, \ldots, w_{n-1}\rangle$ and $L_{v}(u)=u$. Hence we can find a $g\in {\SL}_{n}(\mathbb{Q})$ such that 
$H_v =g^{-1} \left[\begin{array}{c|c}
SO( Q_{v \mid L_v =0})& 0\\\hline
0 & 1
\end{array}\right]g.$

\noindent  Let us define the hypersurface defined over $\mathbb{Q}$ by  $X_{0}=\{Q_0 = L_{0}=0\}= \{ Q_{0 \mid L_0 =0} =0 \}$ and set
$$H_{0}:=g^{-1} \left[\begin{array}{c|c}
SO( Q_{0 \mid L_0 =0})& 0\\\hline
0 & 1
\end{array}\right]g.$$
Clearly we have $X_v = X_0$ and $H_v = H_0$ where $v$ is as chosen above. The form  $Q_{0  | L_{0} =0}$ is nondegenerate thus $H_0$ is semisimple,  it  is also isotropic which implies that  $H_{0}$ is noncompact and isotropic at $v$. Moreover, from the rationality of $Q_{0  | L_{0} =0}$  we  infer that $H_{0} $ is defined over $\mathbb{Q}$. Therefore the triple $(X_0, Q_{0  | L_{0} =0}, H_{0})$ satisfies all the conditions of Theorem \ref{main} unless the connectedness for $H_0$ which is not ensured at all. To remedy to this situation we can consider the connected component of the identity of $H_0$ which we denote $H_{0}^{+}$. The point is that $H_{0}^{+}$ is now connected but it is not obvious that the other properties (isotropy and rationality) are preserved by performing this operation. The key fact is that $H_{0}^{+}$ has finite index in $H_{0}$ thus since $H_{0}^{+}$ is still noncompact and therefore isotropic, in addition $H_{0}^{+}$ is defined over $\mathbb{Q}$ since $H_0$ is (see e.g. Prop. 1.2(b), \cite{borel} and also 2.3.2 \cite{Mbook} and remark 2 just after). A last remark concerns the conservation by the connected component under the central isogeny. This point is quite crucial since at some stage in the proof of the Theorem \ref{main} we need to pass to the universal covering. The fact that $\pi : \widetilde{H_{0}} \rightarrow H_{0}$ is an isogeny defined over $\mathbb{Q}$, we have that $H_{0}^{+}=\pi( \widetilde{H_{0}})^{+} = \pi( \widetilde{H_{0}}^{+})$ (see e.g. Cor. 1.4 (b), \cite{borel}). Hence $\pi$ induces an central isogeny $\pi^{+} : \widetilde{H_{0}}^{+} \rightarrow H_{0}^{+}$ where $ \widetilde{H_{0}}^{+}$ is simply connected. In fact, the latter subgroup is explicitely given by 
$$\widetilde{H_{0}}^{+} =g^{-1} \left[\begin{array}{c|c}
\mathrm{Spin}( Q_{0 \mid L_0 =0})& 0\\\hline
0 & 1
\end{array}\right]g.$$
To sum up, given an arbitrary $\varepsilon >0$,  the triple $(X_0, Q_{0  | L_{0} =0}, H_{0}^{+})$ satisfies all the conditions of Theorem \ref{main}  for $S_1 = \{ v\}$ thus there exists a nonzero  $x\in \mathbb{Z}_{S}^{n}$ such that $x\in X_{S}^{\varepsilon}$  i.e. 
 \begin{center}
 $ |Q_{s}(x)|_s <  \varepsilon $  ~and~ $ |L_{s}(x)|_{s} <  \varepsilon $ ~~ for each $s \in S$. \begin{flushright}
  $\square$
 \end{flushright}
 \end{center}

\section{Proof of Corollary \ref{cor-pairs-2}}

We proceed in the same way as the previous corollary. Let us consider a $Q$ quadratic form and a linear map $M=(L_{1}, \ldots, L_{r})$ satisfying all the assumptions of the Corollary \ref{cor-pairs-2}. For each $s\in S$, we define the following variety $X_{s}=V(Q_{s}, L_{1,s}, \ldots, L_{r,s})$ defined over $\mathbb{Q}_{s}^{n}$, for pratical reasons it is more suitable to see this variety as  $X_s = V(Q_{s \vert M_s=0})$. Geometrically the locus of $X_{s}$ is defined by a quadratic $\{Q_s =0\}$ cutted out by the intersection of the hyperplanes of equations  $\{L_{i,s} =0\}$ $(1\leqslant i\leqslant r)$, the later intersection is just the kernel of $M_{s}$ i.e.  $X_s =V(Q_{s \vert \ker M_s})$. It is assumed that ${\rank} ~M_s = r$ for every $s\in S$, thus $\dim \ker M_s = n-r$, that is to say, $L_{1,s}, \ldots, L_{r,s}$ are linearly independent over $\mathbb{Q}_{s}$. If we denote by $B_s$ the symmetric bilinear form associated to $Q_s$, we can consider the following orthogonal decompositon of the quadratic space $(\mathbb{Q}_{s}^{n}, Q_s)$  with respect to $B$:
$$ (\mathbb{Q}_{s}^{n}, Q_s) = (\ker M_s, Q_{s | \ker M_s}) \bigoplus ((\ker M_s)^{\perp}, Q_{s | (\ker M_s)^{\perp}})$$
and let $\{e_1, \ldots, e_n\}$ a basis of $\mathbb{Q}_{s}^{n}$ adapted to this decompostion above, that is, $M_{s}=\langle e_{1}, \ldots, e_{n-r}\rangle$ and $M_{s}^{\perp}=\langle e_{n-r+1}, \ldots, e_{n}\rangle$ such that $B(e_{k},e_{l})=0$ for all $1\leqslant k \leqslant n-r$ and $n-r+1\leqslant k \leqslant n$. Now for any $s\in S$, we set the following subgroup of ${\SL}_{n}(\mathbb{Q}_s)$  given in the adapted basis by 
$$H_s =  \left[\begin{array}{c|c}
SO( Q_{s \mid \ker M_s})& 0 \\ \hline

0 &~~ I_{r}
\end{array}\right].$$

\noindent We claim that $H_{s}$ leaves invariant both $Q_s$ and $M_{s}$ for each $s\in S$. Indeed let $h\in H_s$, and  $A $ an element of $  SO( Q_{s \mid \ker M_s})$ such that
$$h= \left[\begin{array}{c|c}
A & 0 \\ \hline

0 &~~ I_{r}
\end{array}\right].$$ 
Note that, as given, the range of $A$ is necessarily within the subspace $\ker M_s$.
Let $x$ be a vector in $ \mathbb{Q}_{s}^{n}$ which decomposes into $x=x_{1} + x_{2}$ with $x_{1}\in \ker M_s$ and $x_{2}\in (\ker M_s )^{\perp}$. Therefore $L_s(hx) = L_s(h(x_{1},x_{2})^{t})= L_s (Ax_{1} \oplus x_2) = L_s(Ax_{1}) + L_s(x_2) = 0 + L_s (x_2) = L_s(x_1) + L_s(x_2) =L_s(x)$, thus $L_s$ is $H_{s}$-invariant. In the other hand, using the same decomposition for $x$ we get 
$$Q_s(hx) = Q(Ax_{1}\oplus x_2) = Q_{s|\ker M_s}(Ax_{1})+Q_s(x_2) = Q_{s}(x_1)+ Q_{s}(x_2)=Q(x).$$
In particular, the claim shows that $H_s$ acts linearly on $X_s =V(Q_{s | \ker M_s})$. Now given any constants $\alpha_1, \ldots, \alpha_r$ in $ \mathbb{Q}_{S}$ not all zero, the form $\alpha_{1,s} L_{1,s}+ \ldots+ \alpha_{r,s} L_{r,s}$ is irrational for every place $s\in S\backslash \{v\} $ and proportional to a rational form for $s=v$. Let us set $f_{0} := Q_v + \alpha_{1,v} L_{1,v} + \ldots + \alpha_{1,v} L_{r,v}$, it is clear that $f_0$ is proportional to a rational form and a suitable choice of constants allows us to assume that $f_0$ has rational coefficients. The crucial fact is that $f_{0 | \ker M_v } = Q_{v \mid \ker M_v}$  is also a quadratic form with rational coefficients since $\ker M_v $ is a $\mathbb{Q}$-subspace. Thus if we put 
$$H_0 =  \left[\begin{array}{c|c}
SO( f_{0\mid \ker M_v})& 0 \\ \hline

0 &~~ I_{r}
\end{array}\right]$$
then $H_{0}=H_v$ is a algebraic subgroup of ${\SL}_{n}(\mathbb{Q}_v)$ which is defined over $\mathbb{Q}$ and which acts on $X_0=X_v$. Moreover, for the same reasons as in the previous corollary,  $H_{0}^{+}$ is a connected semisimple algebraic subgroup which is isotropic at $v$ since $Q_{v | \ker M_v}$ is isotropic, in particular $H_{0}^{+}$ has no compact factors. We obtain a rational triple by taking $(X_0, f_{0 | \ker M_v}, H_{0}^{+})$, since it satisfies all the conditions of Theorem \ref{main} for $S_{1}=\{v\}$, we infer that for any $\varepsilon >0$, we can find a nonzero $x\in \mathbb{Z}_{S}^{n}$ such that 
 \begin{center}
$|Q_{s}(x)|_{s} < \varepsilon$ and $ |L_{i,s}(x)|_{s} <  \varepsilon $ ~~ for each $s \in S$ and $1\leqslant i \leqslant r$.
 \end{center} 
To conclude, one has to remark that since $Q_{s}$ is rational and $x\in \mathbb{Z}_{S}^{n}=(\mathbb{A}_{S}\cap \mathbb{Q})^{n}$, $Q_{s}(x)\in \mathbb{Q}$ for every $s\in S$. Since $\mathbb{Q}$ is discrete in every completion, we deduce that if $\varepsilon$ is small enough we can find $x\in \mathbb{Z}_{S}^{n}-\{0\}$ such that 
 \begin{center}
$Q_{s}(x) =0$ and $ |L_{i,s}(x)|_{s} <  \varepsilon $ ~~ for each $s \in S$ and $1\leqslant i \leqslant r$. \ \ \ $\square$
 \end{center}

\section{Existence of rational triples $(X_{0},H_{0},f_{0})$ for general varieties}

In this section we adress some remarks concerning the class of varieties $X$ which falls into the conditions of the theorem. The varieties involved in the theorem are called \textit{complete intersection} in the litterature and they have been subject to extensive research until now and still many problems remains open concerning those projective varieties. From our point of view, we are more concerned with the invariant theory of the space of homogeneous polynomials. \\

\noindent In order to apply the main theorem one has to find a rational  triple  $(X_0, H_0, f_0)$ such that $X_S$ (resp. $H_S$) can be splited in the form $X_{S}= V(f_0) \times X_{S_2}$ (resp. $H_S =H_0 \times H_{S_{2}}$) without loss in generality we assume that $S_{1}=\{v\}$ for some $v\in S_f$ thus $S_2=S\backslash \{v\}$. In particular $H_0$ acts rationally on the $\mathbb{Q}$-hypersurface $X_0=V(f_{0})$.  \\

\subsubsection*{Bounds on the degrees $(d_1, \ldots, d_r)$ of the generators of the ideal of $X$}~\\

$\bullet$  A first constrain is the equality $X_{v}=X_{0}$, that is, in algebraic terms $$V(f_{1,v}, \ldots, f_{r,v}) =V(f_{0}).$$

\noindent Applying Hilbert's Nullstellensatz in an algebraic closure of $\mathbb{Q}_{v}$ yields 
$$ \sqrt{(f_0)} =  \sqrt{(f_{1,v}, \ldots, f_{r,v})}$$
where $\sqrt{J}$ is the radical of an ideal $J$ in $\overline{\mathbb{Q}_{v}}[x_1, \ldots, x_n]$, in particular there exists an integer $\rho>0$ and (homogeneous) polynomials $P_1, \ldots, P_r$ over $\mathbb{Q}_v$ such that 
\begin{equation} \label{nullst}
f_{0}^{\rho} = P_1 f_{1,v} + \ldots + P_{r} f_{r,v}
\end{equation}

\noindent Let us denote by $N_i$ (resp. $d_i$) the total homogeneous of $P_{i}$ (resp. $f_{i,v}$) for each $1 \leqslant i \leqslant r$.  Thus the previous equality reads in terms of degrees as 
\begin{equation}\label{degf0}
\rho \deg f_{0} = \max_{1 \leqslant i \leqslant r} \{ N_{i} + d_{i} \}.
\end{equation}

\noindent Let us assume that the degrees are ordered as follows $N_{1} \geqslant N_{2} \geqslant \ldots \geq N_{r}$ and $d_{1} \geqslant d_{2} \geqslant \ldots \geqslant d_{r}$, then (\ref{degf0}) reads
\begin{equation}
\rho \deg f_{0} =  N_{1} + d_{1}.
\end{equation}

\noindent When $d_i \neq 2$ $(1\leqslant i \leqslant r)$, upper bounds for $\rho$ can be effectively computed, the following sharp estimates for $\rho$ are due to J. Koll\'{a}r (Corollary 1.7, \cite{kollar}.)

\begin{equation} \label{rho}
 \rho \leq \left\lbrace 
\begin{array}{ccc}
d_{1}d_{2} \ldots d_r & \mathrm{if} & r \leq n \\
 d_{1}d_{2} \ldots d_{n-1} d_r & \mathrm{if} & r >n
\end{array}\right.
\end{equation}
and for each $1 \leqslant i \leqslant r$
\begin{equation}\label{NplusP}
N_i + d_{i} = \deg(P_{i}f_{i,v}) \leq \left\lbrace 
\begin{array}{ccc}
(1+d_0 )~d_{1}d_{2} \ldots d_r & \mathrm{if} & r \leq n \\
~~~(1+d_0 )~ d_{1}d_{2} \ldots d_{n-1} d_r & \mathrm{if} & r >n
\end{array}\right.
\end{equation}
 where $d_0$ denotes $\deg f_0$. To sum up,  saying that $X_{v}$ equals the hypersurface $X_{0}=V(f_{0})$ amounts to find some polynomials $(P_{j})_{1 \leqslant i \leqslant r}$ such that 
 $$f_{0}^{\rho} = P_1 f_{1,v} + \ldots + P_{r} f_{r,v}$$

\noindent  where $\rho$ satisfies condition (\ref{rho}) and after (\ref{NplusP}) we get that for every $1 \leqslant i \leqslant r$, 
\begin{equation}
0 \leq \deg(P_{j}) \leq  \left\lbrace 
\begin{array}{ccc}
(1+d_0 )~d_{1}d_{2} \ldots d_r -d_{j}& \mathrm{if} & r \leq n \\
~~~(1+d_0 )~ d_{1}d_{2} \ldots d_{n-1} d_r  -d_{j}& \mathrm{if} & r >n.
\end{array}\right.
\end{equation}

\noindent \textit{The equidimensional case}:  If we assume that $d_1 = \ldots =d_r =d \neq 2$, we obtain the following bounds for every $1 \leqslant i \leqslant r$, 
\begin{equation}\label{equidim}
0 \leq \deg(P_{j}) \leq  \left\lbrace 
\begin{array}{ccc}
(1+d_0 )~d^r -d& \mathrm{if} & r \leq n \\
~~~(1+d_0 )~ d^{n} -d& \mathrm{if} & r >n.
\end{array}\right.
\end{equation}

\noindent If we want that the right hand in (\ref{degf0}), that is, $f_{0}^{\rho}$ to have only degree $d$, then necessarily  $\rho=1$ and from (\ref{equidim}) we get that the polynomials $P_i$ should be constant polynomials $\alpha_{i}$ $(1\leqslant i \leqslant r)$

\begin{equation}\label{f0-lin}
f_{0} = \alpha_1 f_{1,v} + \ldots + \alpha_r f_{r,v}. 
\end{equation}
~\\
\noindent $\bullet$ Rationality conditions $X_v = X_{0}=V(f_{0})$ with $f_0$ rational. The relation (\ref{nullst}) shows  that if the polynomial $f_0$ has rational coefficients then some linear combination  (over $\mathbb{Q}_v$) of the $f_{i,v}$'s must be rational. In particular,  this and (\ref{f0-lin}) explain why the condition (3) in both corollaries \ref{cor-pairs} and \ref{cor-pairs-2} is necessary. \\

\subsubsection*{From complete intersections towards invariant group of the ideal of definition}
~  \\

\noindent  The main issue is that given a variety $X$ defined over a field of characteristic zero $K$, say a complete intersection, to find the largest subgroup of $G={\SL}_{n | K}$ which acts trivially on the ideal of definition $I_X$. For instance, let us be given a complete intersection $X=V(f_1, \ldots, f_r)$ where $f_1, \ldots, f_r$ are homogeneous polynomial of degrees $d_{1} \leq \ldots \leq d_r$. The ideal of definition of $X$ is given by $I_X = \langle f_1, \ldots, f_r \rangle$. The central role is played by the stabilizer $H$ of the ideal which is defined to be 
$$H_{s} =\bigcap_{1\leqslant i \leqslant r} \{ g\in G ~|~  g.f_{i} =f_{i}\}.$$  \\
The ideal stabilizer $H_s$ obviously does act on $X_s$ for every $s\in S$, and it gives an action of $H_{S}$ on $X_S$ induces by the usual of $G=\mathrm{PSL_{n}}$ on the vector space of homogeneous polynomials in $\mathbb{Q}_{s}[x_1, \ldots, x_n]$. In particular $H_s$ contains $\mathrm{Aut}_{G}(X)$ the group automorphism of $X$ under $G$ which is the pointwise stablizer of $X_{s}$ under $G$.\\
The ideal stabilizer is an algebraic subgroup of $G$ given by the following equations in the variables $(g_{i,j})$
\begin{center}
$ f_{k}(gx)=f_k(x)  $  and  $\det (g_{i,j}) -1 = 0$.
\end{center}

\noindent Let us try to solve those equation with $g =(g_{i,j})_{i,j}$, for this let us explicit the coefficients of the $f_k$ and assume that they are of form
$$  f_k(x_1,\ldots,x_{n}) = \sum_{|\alpha|=d_k} a_{\alpha}^{(k)} x_{1}^{\alpha_1}\ldots x_{n}^{\alpha_n}.$$
Therefore we have $n $ equations which takes places in $K[x_1,\ldots,x_n]_{(d_k)}[(g_{i,j})_{i,j}]$
$$ \sum_{|\alpha|=d_k} a_{\alpha}^{(k)} \left((\sum_{j_1=1}^{n} g_{1j_1}x_{j_1})^{\alpha_1}\ldots (\sum_{j_n=1}^{n} g_{nj_n}x_{j_n})^{\alpha_n}-x_{1}^{\alpha_1}\ldots x_{n}^{\alpha_n}\right) =0.$$

$$ \sum_{|\alpha|=d_k} a_{\alpha}^{(k)} (\sum_{|\beta|=\alpha_1} \left(\begin{array}{c}
\alpha_1 \\ \beta
\end{array}\right) (g_{11}x_{1})^{\beta_1}   \ldots(g_{1n} x_{n})^{\beta_{n}}) \ldots (\sum_{|\beta|=\alpha_n} \left(\begin{array}{c}
\alpha_n \\ \beta
\end{array}\right) (g_{n1}x_{1})^{\beta_1}   \ldots(g_{nn} x_{n})^{\beta_{n}}   ) $$
$$ = \sum_{|\alpha|=d_k} a_{\alpha}^{(k)} x_{1}^{\alpha_1}\ldots x_{n}^{\alpha_n}.$$

We do not need to go further to observe that such compuations unless we are dealing with  low degrees (i.e. $d=2,3$) leads to tremendous compuations and trying to obtain the ideal stabilizer in such a way is quite compromized. Using elimination when $d=2,3$ one could provide  the required invariant groups where the solutions $g=(g_{i,j})$ are given by functions of the coefficients $a_{\alpha}^{(k)}$ of the $f_k$'s. A last step would be to decide if the invariant group is compact/semisimple/isotropic which could ask some additional efforts, this task is  crucial in order to apply our main theorem.\\

\noindent  \textit{Final Comments}\\
The main problem is to determine the ideal stabilizer $H$ of a given projective variety $X=V(f_1, \ldots,f_r)$. This consists in finding a subgroup $H$ such that $$K[x_1, \ldots, x_n]^{H} = K[f_1, \ldots, f_r].$$ 
This question is dual to the \textit{Invariant theory}, indeed in invariant theory we fix a group $G$ and we try to understand the ring of invariants $k[X]^{G}$ of the variety $X$. This theory has reached a good level of maturity, notably with the rise of the \textit{geometric invariant theory }(G.I.T. \cite{mumford}) and more recently with the theory of \textit{prehomogeneous spaces} for which we hope that we could derive an analog of the work of A. Yukie (\cite{yukie}) using our main theorem. In terms of category, the invariant theory is an attempt to understand the image of the functor 
$$ \begin{array}{cccc}
F(X):  & \mathbf{Grps} & \rightarrow & \mathbf{Rings}\\
    & G & \rightarrow & k[X]^{G}.
\end{array}
$$

The so called \textit{Inverse Invariant theory} consists the dual situation, namely the image of the following functor given a fixed projective variety

$$ \begin{array}{cccc}
F(X)^{\ast}:  &\mathbf{Rings} & \rightarrow & \mathbf{Grps} \\
    & A & \rightarrow & G
\end{array}
$$
where we define the functor $F(X)^{\ast}$ as follows $F(X)^{\ast}(A)=G$ if $A=  k[X]^{G}$. As far as we know, this theory has been only  developped for finite groups and in particular for linear groups over \textit{finite fields}. The latter has been studied in detail by Neusel using tools from algebraic topology such as Steenrod operations. It should be very interesting to have such a theory for linear groups over fields in null characteristic and to have a criterion which ensures that the group obtained is reductive or/and noncompact. If we could have such theory it would open a large range of applications and more particularly for solving the diophantine inequalities in (\ref{ineq}).

\end{document}